\documentclass[11pt]{article}

\title{More on the Bernoulli*-Taylor formula for extended umbral calculus}
\author{A.K.Kwa\'sniewski\\  
\\ High School of Mathematics and Applied Informatics\\
PL - 15-021 Bialystok , ul.Kamienna 17,  Poland
\\e-mail: kwandr@wp.pl}

\usepackage{amsmath,amsthm}

\chardef\bslash=`\\ 
\begin{document}
\maketitle

\begin{abstract}
\noindent One presents here the $*_{\psi}$-Bernoulli-Taylor*
formula of a new sort  with the rest term of  the Cauchy type
recently derived by the author in the case of the so called
$\psi$-difference calculus which constitutes the representative
for the purpose case of extended umbral calculus. The central
importance of such a type formulas is beyond any doubt - and
recent publications do confirm this historically established
experience.Its links via umbrality to combinatorics are known at
least since Rota and Mullin-Rota source papers then up to recently
extended by many authors to be indicated in the sequel.

\end{abstract}

\vspace{2mm}


\vspace{2mm}
AMS S.C.   (2000)   05A40, 81S05,  01A45, 01A50, 01A61
\vspace{2mm}

* see below : a historical remark based on Academician N.Y.Sonin article published in Petersburg in 19-th century.
We owe this information and the article to Professor O.V.Viscov from Moscow.\\

\vspace{2mm}

** see: C.~Graves,  {\it On the principles which regulate the
interchange of symbols in certain symbolic equations}, Proc. Royal
Irish Academy {\bf 6}(1853--1857), 144-152

\section{One Historical Remark} Here are the famous examples of expansion
$$\partial_0=\sum_{n=1}^{\infty}\frac{x^{n-1}}{n!}\frac{d^n}{dx^n}$$
or  $$\epsilon_0=\sum_{n=0}^{\infty}(-1)^n
\frac{x^{n}}{n!}\frac{d^n}{dx^n}$$ where $\partial_0$ is the
divided difference operator while $\epsilon_0$ is at the zero
point evaluation functional. If one compares these with
\textit{"series universalissima"} of J.Bernoulli from \textit{Acta
Erudicorum} (1694) (see commentaries in \cite{12}) and with
$$exp\{yD\}=\sum_{k=0}^{\infty}\frac{y^kD^k}{k!},\ \
D=\frac{d}{dx},$$ then confrontation with B.Taylor's {\em
"Methodus incrementorum directa et inversa"} (1715), London;
entitles one to call the expansion formulas considered in this
note {\em "Bernoulli - Taylor formulas"} or (for $n \rightarrow
\infty$) {\em "Bernoulli - Taylor series" }\cite{1}.

\vspace{2mm}

\textbf{Information:}\quad Johann Bernoulli was elected a fellow
of the academy of St Petersburg.Johann Bernoulli - the Discoverer
of \textit{Series Universalissima} was "\textbf{Archimedes of his
age}" and this is indeed inscribed on his tombstone.

\vspace{2mm}

\section{Introduction}

\textbf{A.} From here now $\psi$ denotes an extension of
$\langle\frac{1}{n!}\rangle_{n\geq 0}$ sequence to quite arbitrary
one (the so called "admissible"- see: Markowsky) and the specific
choices are for example : Fibonomialy  -extended ($\langle F_n
\rangle$ - Fibonacci sequence )
$\langle\frac{1}{F_n!}\rangle_{n\geq 0}$   or just "the usual"
$\langle\frac{1}{n!}\rangle_{n\geq 0}$ or Gauss $q$-extended
$\langle\frac{1}{n_q!}\rangle_{n\geq 0}$ admissible sequences of
extended umbral operator calculus - see more in
\cite{16,13,14,15}.

The simplicity of the first steps to be done while identifying
general properties of such [6,7,8,16,13-15]  $\psi$-extensions
consists in natural notation i.e. here - in writing objects of
these extensions in mnemonic convenient \textbf{upside down
notation} \cite{15} , \cite{14}

\begin{equation}\label{eq31}
\frac {\psi_{(n-1)}}{\psi_n}\equiv n_\psi,
n_\psi!=n_\psi(n-1)_\psi!, n>0 ,   x_{\psi}\equiv \frac
{\psi{(x-1)}}{\psi(x)} ,
\end{equation}

\begin{equation}\label{eq32}
x_{\psi}^{\underline{k}}=x_{\psi}(x-1)_\psi(x-2)_{\psi}...(x-k+1)_{\psi}
\end{equation}

\begin{equation}\label{eq33}
x_{\psi}(x-1)_{\psi}...(x-k+1)_{\psi}=
\frac{\psi(x-1)\psi(x-2)...\psi(x-k)} {\psi(x)
\psi(x-1)...\psi(x-k +1)} .
\end{equation}

If one writes the above in the form $x_{\psi} \equiv \frac
{\psi{(x-1)}}{\psi(x)}\equiv \Phi(x)\equiv\Phi_x\equiv x_{\Phi}$ ,
one sees that the name upside down notation is legitimate.

 You may consult \cite{10} , \cite{9} ,[14,15] for further
development and use of this notation.

With such an extension we frequently though not always may
"$\psi$-mnemonic" repeat with exactly the same simplicity and
beauty most of what was done by Rota \cite{10,15}. Accordingly the
extension of notions and formulas with its elementary essential
content and context to general case of $\psi$- umbral instead of
umbral or $q$-umbral calculi case only - is sometimes automatic
\cite{10} , \cite{14,15} (see corresponding earlier references
there and necessary definitions).

\textbf{B.} While deriving  the Bernoulli-Taylor $\psi$-formula
one is tempted to adapt  the ingenious  Viskov`s method  [2] of
arriving to formulas of such type  for various pairs of
operations. In our case these would be $\psi$-differentiation and
$\psi$-integration (see: Appendix). However  straightforward
application of  Viscov methods in  $\psi$-extensions of umbral
calculus leads to sequences which are not normal (Ward) hence a
new invention is needed. This expected and verified here invention
is the new specific $*_{\psi}$ product of analytic functions or
formal series. This note is based on [3]  where the  derivation of
this new form of Bernoulli-Taylor $*_{\psi}$ - formula was
delivered due to the use of a specific $*_{\psi}$ product of
formal series.

\section{Classical  Bernoulli-Taylor formulas with the rest term of  the Cauchy type by Viskov method}
Let us consider the obvious identity
\begin{equation} \label{eq1}
  \sum_{k=0}^{n} (\alpha_{k}-\alpha_{k+1})=\alpha_{0}-\alpha_{n+1}
\end{equation}
in which (\ref{eq1}) we now put
$\alpha_{k}=a^k b^k ; a,b\in \mathcal{A}.$
$\mathcal{A}$ is an associative algebra with unity over the field F=R,C.
Then we get
\begin{equation} \label{eq2}
\sum_{k=0}^n a^k (1-ab) b^k = 1-a^{n+1}b^{
n+1} ; a,b \in \mathcal{A}\
\end{equation}
Numerous choices of  $a,b \in \mathcal{A}$ result in many important
specifications of (\ref{eq2})

{\bf Example 1.} Let $\mathcal{F}$ denotes the linear space of sufficiently
smooth functions ${\it f: F \longrightarrow F}$. Let
\begin{eqnarray} \label{eq3}
a: \mathcal{F}\longrightarrow \mathcal{F};\quad &&(af)(x)= \int_{a}^{b}
f(t)dt,\nonumber \\
b: \mathcal{F}\longrightarrow \mathcal{F};\quad  &&(bf)(x) =
(\frac{d}{dx}f)(x);\\ l:\mathcal{F}\longrightarrow \mathcal{F};
\quad&& (lf)(x)=f(x)\nonumber.\end{eqnarray}
Then [b,a]=1-ab=$\varepsilon_\alpha$ where $\varepsilon_\alpha$ is evaluation
functional on $\mathcal{F}$ i.e.
\begin{equation}
\varepsilon_\alpha(f)=f(\alpha)
\end{equation}
\\
Using now the text-book integral Cauchy formula ($k>0$)
\begin{equation}
(a^k f)(x)= \int_a^x\frac{(x-t)^{k-1}}{(k-1)!} f(t)dt,\
\end{equation}
\\
and under the choice (\ref{eq3}) one gets from (\ref{eq2}) the well-known
Bernoulli-Taylor formula
\begin{equation} \label{eq6}
f(x)= \sum_{k=0}^n\frac{(x-\alpha)^k}{k!} f^{(k)} (\alpha)+R_{n+1}(x)
\end{equation}
\\
with the rest term $R_{n+1}(x)$ in the Cauchy form
\begin{equation} \label{eq7}
R_{n+1}(x)=\int_a^x\frac{(x-t)^{n}}{n!} f^{(n+1)}(t)dt
\end{equation}
\\
{\bf Example 2.} \cite{1} Let $\mathcal{F}$ denotes the linear space of
functions $f:Z_{+}\longrightarrow F; Z_{+}=N\cup \{0\}$. Let

\begin{eqnarray} \label{eq8}
a: Z_{+}\longrightarrow \mathcal{F};\quad &&(af)(x)= \sum_{k=0}^{x-1} f(k),
\nonumber \\
b: Z_{+}\longrightarrow \mathcal{F};\quad  &&(bf)(x) = f(x+1)-f(x),\\
l:Z_{+}\longrightarrow \mathcal{F};\quad&& (lf)(x)=f(x)\nonumber.
\end{eqnarray}
It is easy to see that  [b,a]=1-ab=$\varepsilon_0$ where
$\varepsilon_0$ is evaluation functional  i.e.
$\varepsilon_0(f)=f(0).$ $b=\Delta$ is the standard difference
operator with its left inverse definite summation operator a. The
corresponding $\Delta$ - calculus Cauchy formula is also known
(see formula (31 p.310 in \cite{5});
\begin{equation}
(a^k f)(x)= \sum_{r=0}^{x-1} \frac{(x-r-1)^{\underline {k-1}}}{(k-1)!}f(r);
k>0
\end{equation}
where $x^{\underline n}=x(x-1)(x-2)...(x-n+1).$
\\
Under the choice (\ref{eq8}) one gets from (\ref{eq2}) the $\Delta$ -
calculus Bernoulli - Taylor fomula \cite{1}
\begin{equation}
f(x)=\sum_{k=0}^n \frac {x^{\underline k}}{k!}(\Delta^k f)(0)+R_{n+1}(x)
\end{equation}
with the rest term $R_{n+1}(x)$ in the Cauchy $\Delta$ form
\begin{equation}
R_{n+1}(x)= \sum_{r=0}^{x-1} \frac{(x-r-1)^{\underline {n}}}{n!}
(\Delta^{n+1}f)(r);
\end{equation}

\section{$"\ast_\psi$ realization" of Bernoulli
identity.}  Now \textit{a specifically new} form of the
Bernoulli-Taylor formula with the rest term of the Cauchy type as
well as Bernoulli-Taylor series is to be supplied in the case of
$\psi$-difference \textit{umbral} calculus (see [5-8] and [9,10]
and references therein). For that to do we use natural
$\psi$-umbral representation [13,14] of Graves-Heisenberg-Weyl
(GHW) algebra  [11,12]  generators $\hat{p}$ and $\hat{q}$ and
then we use Bernoulli identity (\ref{eq12})
\begin{equation}  \label{eq12}
\hat p \sum_{k=0}^n \frac {(-\hat q)^k \hat p^k}{k!}=\frac{(-\hat q)^n \hat
p^{n+1}}{n!}
\end{equation}
derived by Viskov from (\ref{eq1}) under the substitution (see
(28) in [2])
$$\alpha_0=0,\ \alpha_k=(-1)^k{(\hat q)^{k-1} \hat p^k}{(k-1)!},\
k=1,2,...$$
due to $\hat p \hat q^n= \hat q^n \hat p +n \hat
q^{n-1}$  (n=1,2,...) resulting by induction from
\begin{equation} \label{eq13}
[\hat p, \hat q]=1
\end{equation}
\\{\bf Example 1.} The choice $\hat p=D \equiv \frac{d}{dx}$ and $\hat q
=\hat x-y, y\in F; \hat x f(x)=xf(x)$ after substitution into
Bernoulli identity (\ref{eq12}) and integration $\int_\alpha^x dt$
gives the Bernoulli - Taylor formula (\ref{eq6}).
\\{\bf Example 2.} The choice \cite{2} $\hat p=\Delta$  and $\hat q=
\hat x\circ E^{-1}$ where $E^{\alpha}f(x)= f(x+\alpha)$ after
substitution into Bernoulli identity (\ref{eq12}) and $"\Delta$ -
integration" $\sum_{r=0}^{\alpha-1}$ gives the Bernoulli - Mac
laurin formula of the following form ($\alpha, x \in {\bf Z},
\bigtriangledown =1-E^{-1})$ with the rest term $R_{n+1}(x)$
\begin{equation}
f(0)=\sum_{k=0}^{n} \frac{\alpha^{\underline{k}}}{k!}(-1)^{k+1}
(\bigtriangledown ^k f)(\alpha)+R_{n+1}(\alpha);
\end{equation}
\begin{equation}
R_{n+1}(\alpha)=(-1)^n\sum_{r=0}^{\alpha-1} \frac{r^{\underline n}}{n!}
(\bigtriangledown^{n+1}f)(r+1).
\end{equation}
\\{\bf Example 3.} Here  $f^{(k)}\equiv \partial_{\psi}^k f$ and $f(x)*_{\psi}g(x)\equiv f(\hat x_{\psi})g(x)$
 - see Appendix. The choice $\hat p =\partial_{\psi}$ and $ \hat q
= \hat z_{\psi}$  $(z=x-y)$ where $\hat x_{\psi}
x^n=\frac{n+1}{(n+1)_{\psi}}x^{n+1}$ after substitution into
Bernoulli identity (\ref{eq12}) and "$\partial_{\psi}$ -
integration" $\int_\alpha^x d_{\psi}t$ (see: Appendix) gives
another Bernoulli - Taylor $\psi$-formula of the form:
\begin{equation}
f(x)=\sum_{k=0}^n
\frac{1}{k!}(x-\alpha)^{k_{\ast_{\psi}}}\ast_{\psi}
f^{(k)}(\alpha)+ R_{n+1}(x)
\end{equation}
\\with the rest term $R_{n+1}(x)$ in the Cauchy-form
\begin{equation}
R_{n+1}(x)=\frac{1}{n!}\int_\alpha^x d_qt (x-t)^{n_{\ast_q}}\ast_q f^{(n+1)}
(t)dt
\end{equation}
\\In the above notation $x^{0*_{\psi}}=1,\ x^{n*_{\psi}}\equiv x*_{\psi} (x^{(n-1)*_{\psi}})=
x*_{\psi} ...*_{\psi} x=\frac{n!}{n_{\psi}!}x^n;\ n\geq 0.$
\\Naturally $\partial_{\psi} x^{n*_{\psi}}=nx^{(n-1)*_{\psi}}$ and in general $f,g$ - may be
formal series for which
\begin{equation}
\partial_{\psi}(f*_{\psi} g)=(Df)*_{\psi} g+f*_{\psi}(\partial_{\psi} g)
\end{equation}
\\i.e. Leibniz $*_{\psi}$ rule holds \cite{13,14,15}.
\\ Summary: These another forms of both the Bernoulli -Taylor
formula with the rest term of the Cauchy type \cite{3} as well as
Bernoulli - Taylor series are quite easily handy due to the
technique developed in \cite{13,14} where one may find more on
$*_\psi$ product devised perfectly suitable for the Ward's {\em
"calculus of sequences"} \cite{6} or more exactly $*_\psi$ is
devised perfectly suitable for the so-called $\psi$ - extension on
Finite Operator Calculus of Rota (see \cite{9,10,14,15} and
references therein)
\section{Appendix}  \textbf{$*_{\psi}$ product}\\
Let $n-{\psi}\equiv\psi_n$; $\psi_n\neq 0$: $n>0$.  Let
$\partial_\psi$ be a linear operator acting on formal series  and
defined accordingly by $\partial_\psi x^n=n_\psi x^{n-1}$. \\We
introduce now a intuition appealing
$\partial_\psi$-difference-ization rules for a specific new
$*_\psi$ product of functions or formal series. This $*_\psi$
product is what we  call:  the $\psi$-multiplication of functions
or formal series as specified below.
\\\textbf{Notation A.1.}\\
$x \ast _{\psi} x^{n} = \hat {x}_{\psi}(x^{n}) = \frac{{\left( {n
+ 1} \right)}}{{\left( {n + 1} \right)_{\psi} } }x^{n + 1};\quad n
\geq 0$ \; hence $x \ast _{\psi} 1 = (1_{\psi})^{-1}  x \not
\equiv x $ therefore $x  \ast_{\psi} \alpha 1 = \alpha 1 \ast
_{\psi}  x = x  \ast _{\psi} \alpha = \alpha \ast _{\psi}  x =
\alpha  (1_{\psi})^{-1} x$ and $\forall x, \alpha\in F$; $f(x)
\ast _{\psi} x^{n} = f(\hat {x}_{\psi})x^{n}$.
\\For $k \ne n $ \; x$^{n}
\ast _{\psi} $ x$^{k} \ne$ x$^{k} \ast _{\psi} $ x$^{n}$ as well
as x$^{n} \ast _{\psi} $ x$^{k} \ne$ x$^{n+ k}$ - in general.
\\In order to facilitate the  formulation  of observations  accounted
for on the basis of $\psi$-calculus representation of  GHW algebra
we shall use what follows. \\

\textbf{Definition A.1.} With Notation A.1. adopted define the
$*_\psi$ powers of $x$ according to
 $x^{n \ast _{\psi} } \equiv $ x $ \ast _{\psi} x^{\left( {n - 1} \right)
\ast _{\psi} } = \hat {x}_{\psi} (x^{\left( {n - 1} \right)\ast
_{\psi }} ) = $ x $ \ast _{\psi} $ x $ \ast _{\psi} $ ... $ \ast
_{\psi} $ x $=\frac{n!}{n_{\psi}  !}x^{n};\quad n \geq 0$. Note
that $x^{n\ast _{\psi} }  \ast _{\psi} x^{k\ast _{\psi} } =
\frac{{n!}}{{n_{\psi}  !}} x^{\left( {n + k} \right)\ast _{\psi} }
\ne x^{k\ast _{\psi} }  \ast _{\psi} x^{n\ast _{\psi} }  =
\frac{{k!}}{{k_{\psi}  !}}x^{\left( {n + k} \right)\ast _{\psi} }
$ for
$k  \ne n$ and $x^{0\ast _{\psi} }=1$.\\

This noncommutative $\psi$-product $ \ast _{\psi}$ is devised so
as to ensure the following observations.\\

\textbf{Observation A.1}
\begin{enumerate}
\renewcommand{\labelenumi}{\em \alph{enumi})}
\item $\partial _{\psi}  x^{n\ast _{\psi} } = n x^{\left( {n - 1}
 \right)\ast _{\psi} } $;\; $n \ge 0$
\item  {\bf exp}$_{\psi}  $[$\alpha ${\it x}] $ \equiv ${\bf exp}
\{$\alpha \hat {x}_{\psi}  $\}{\bf 1}
\item {\bf exp} [$\alpha x$] $ \ast _{\psi} $
({\bf exp}$_{\psi} $\{$\beta \hat {x}_{\psi}  $\}{\bf 1}) = ({\bf
exp}$_{\psi}  $\{[$\alpha +\beta $]$\hat {x}_{\psi}  $\}){\bf 1}
\item
  $\partial _{\psi} (x^{k} \ast _{\psi} \quad x^{n\ast _{\psi} } ) =
(D x^{k}) \ast _{\psi} x^{n \ast _{\psi} } + x^{k} \ast _{\psi}
(\partial _{\psi}  x^{n\ast _{\psi} })$
\item \label{e} $\partial
_{\psi} ( f \ast _{\psi} g) = ( Df) \ast _{\psi} g  + f \ast
_{\psi} (\partial _{\psi} g)$ ; $f,g$  - formal series \item
\label{f} $f( \hat {x}_{\psi}) g (\hat {x}_{\psi} )$ {\bf 1} $=
f(x) \ast _{\psi} \tilde {g} (x)$ ; $\tilde {g} (x) = g(\hat
{x}_{\psi})${\bf 1}.
\end{enumerate}

\textbf{Umbral "$\sim$" Note:} $\tilde{g}(x)=g(\hat
{x}_{\psi}){\bf 1}$ defines the map $\sim  : g \mapsto \tilde{g}$
i.e. $ \sim : P\mapsto P $ which is an umbral operator.

\vspace{2mm}

\textbf{$\psi$-Integration} Let: $\partial_o x^n=x^{n-1}$. The
linear operator $\partial_o$ is identical with divided difference
operator. Let $\hat {Q}f(x)f(qx)$. Recall also that to the
"$\partial _{q} $ difference-ization" there corresponds the
$q$-integration which is a right inverse operation to
"$q$-difference-ization". Namely

\begin{equation}\label{eqa2}
 F\left( {z} \right): \equiv \left( {\int_{q} \varphi
}  \right)\left( {z} \right): = \left( {1 - q}
\right)z\sum\limits_{k = 0}^{\infty}  {\varphi \left( {q^{k}z}
\right)q^{k}}
\end{equation}
i.e.

\begin{multline}\label{eqa3}
 F\left( {z} \right) \equiv \left( {\int_{q} \varphi }
\right)\left( {z} \right) = \left( {1 - q} \right)z\left(
{\sum\limits_{k =
0}^{\infty}  {q^{k}\hat {Q}^{k}\varphi} }  \right)\left( {z} \right) =\\
=\left( {\left( {1 - q} \right)z\frac{{1}}{{1 - q\hat
{Q}}}\varphi} \right)\left( {z} \right).
\end{multline}
Of course
\begin{equation}\label{eqa4}
\partial _{q} \circ \int_{q} = id
\end{equation}
as
\begin{equation}\label{eqa5}
\frac{{1 - q\hat{Q}}}{{\left( {1 - q} \right)}}\partial _{0}
\left( {\left( {1 - q}\right)\hat {z}\frac{{1}}{{1 - q\hat {Q}}}}
\right)=id.
\end{equation}
Naturally (\ref{eqa5}) might serve to define a right inverse
operation to "$q$-difference-ization"
 $\left( {\partial _{q} \varphi}  \right)\left( {x} \right) = \frac{{1 - q\hat
{Q}}}{{\left( {1 - q} \right)}}\partial _{0} \varphi \left( {x}
\right)$ and consequently  the  "$q$-integration" as represented
by (\ref{eqa2}) and (\ref{eqa3}).  As it is well known  the
definite $q$-integral is an numerical approximation of  the
definite integral obtained in the $q \to 1$ limit.
\\Finally we introduce the analogous representation for $\partial_\psi$ difference-ization

\begin{equation}\label{eqa6}
  \partial_\psi=\hat n_\psi \partial_o;\ \hat n_\psi
  x^{n-1}=n_\psi x^{n-1};\ n\ge 1
\end{equation}
Then
\begin{equation}\label{eqa7}
  \int_\psi x^n=\left(\hat x \frac{1}{\hat
  n_\psi}\right) x^n=\frac{1}{(n+1)_\psi}x^{n+1};\ n\ge 0
\end{equation}
and of course  $\left(\int_\psi \equiv \int d_\psi \right)$
\begin{equation}\label{eqa8}
  \partial_\psi \circ \int_\psi=id
\end{equation}
Naturally $$ \partial_\psi \circ \int_a^x f(t)d_\psi t=f(x)$$ The
formula of   "per partes"  $\psi$-integration  is easily
obtainable from (Observation A.1 e) and it  reads:
\begin{equation}\label{eqa9}
  \int_a^b(f*_\psi \partial_\psi g)(t)d_\psi t=[(f*_\psi
  g)(t)]_a^b - \int_a^b(Df*_\psi g)(t)d_\psi t
\end{equation}
\\
\textbf{Closing Remarks:}
\\
\textbf{I.}  All these above may be quite easily extended \cite
{15} to the case of any $Q\in End(P)$ linear operator that reduces
by one the degree of each polynomial \cite {16}. Namely one
introduces \cite {15}:

\textbf{Definition A.2.}
    $$\hat {x}_{Q} \in End(P), \hat x_{Q}: F[x] \to F[x] $$
     such that  $(x^{n}) = \frac{{\left( {n
+ 1} \right)}}{{\left( {n + 1} \right)_{\psi} } }q_{n + 1}; n \geq
0;$ where $Qq_n=nq_{n-1}$.
Then   $\star_Q$ product of formal
series  and  $Q$-integration are defined almost mnemonic
analogously.

\textbf{II.} In  1937   Jean  Delsarte  [17]  had derived the
general Bernoulli-Taylor formula for a  class of linear operators
$\delta$ including  linear operators that reduce by one the degree
of each polynomial. The rest term of the Cauchy-like type in his
Taylor formula (I) is given in terms of the unique solution of a
first order partial differential equation in two real variables.
This first order partial differential equation  is determined by
the choice of  the linear operator $\delta$  and the function  f
under expansion. In  our   Bernoulli-Taylor -formula (16)-(17) or
in its straightforward $\star_Q$     product of formal series and
$Q$-integration generalization -  there is no need to solve any
partial differential equation.

\textbf{III.}   In  [18] (1941) Professor  J. F. Steffensen - the
Master of polynomials application to actuarial problems
$$(see:http://www.math.ku.dk/arkivet/jfsteff/stfarkiv.htm)$$
supplied a remarcable derivation of another Bernoulli-Taylor
formula with the rest of "Q-Cauchy type" in the example presenting
the "Abel poweroids"

\textbf{IV.} The recent paper [19] (2003) by Mourad E. H. Ismail,
Denis Stanton may serve as a kind of indication  for pursuing
further investigation. There the authors have established two new
q-analogues of a Taylor series expansion for polynomials using
special Askey-Wilson polynomial bases. As "byproducts" their
important paper  includes also  new summation theorems,quadratic
transformations for q-series and new results on a q-exponential
function.

\textbf{V.}  Let us also draw an attention to two more different
publication on the subject which are the ones referred to as
[20,21].The $q$-Bernoulli theorems are named here and there above
as $q$-Taylor theorems. The corresponding $(q,h)$-Bernoulli
theorem for the $\partial _{q,h}$-difference calculus of Hahn [22]
might be also obtained as the the one $(q,0)$-Bernoulli i.e.
$q$-Bernoulli theorem constituting here the special case the
Viskov method [2] application. This is so because the $\partial
_{q,h}$-difference calculus of Hahn [22] may be reduced to
$q$-calculus of {\bf Thomae-Jackson} [5,23] due to the following
observation . Let
$$ h \in R ,(E_{q,h}\varphi)(x)=\varphi (qx+h)$$
and  let
\begin{equation}\label{eqa10}
(\partial _{q,h}\varphi)(x)= \frac{\varphi(x)-
\varphi(qx+h)}{(1-q)x-h}
\end{equation}

\vspace{2mm}

Then (see Hann [5]and [22])
\begin{equation}\label{eq11}
\partial_{q,h} = E_{1,{\frac{{-h}}{{1-q}}}} \partial_q E_{1,{\frac
{h}{{1-q}}}}.
\end{equation}
Due to (30) it is easy now to derive corresponding formulas
including Bernoulli-Taylor $\partial _{q,h}$-formula obtained  in
[24] by the Viskov method [2]  which for
$$q\rightarrow 1 , h\rightarrow 0 $$
recovers the content of one of the examples in [2] , while for
$$q \rightarrow 1 , h \rightarrow 1$$
one recovers the content of the another example in [2]. The case
$h \rightarrow 0 $ is included in the formulas of $q$-calculus of
Thomae-Jackson easy to be specified : see  [22] (see also
thousands of up-date references there). For Bernoulli- Taylor
Formula (presented during $PTM$ - Convention Lodz - $2002$) :
contact [23] for its recent version.

\vspace{1mm}

The comparison of the all above quoted ways to arrive at extended
Bernoulli formulae  we leave for another exhibition of similar
investigations.

\textbf{VI.} As indicated right after \textbf{Observation A.1.e)}
the rule $\tilde{g}(x)=g(\hat {x}_{\psi}){\bf 1}$ defines the map
$\sim  : g \mapsto \tilde{g}$ which is an umbral operator   $ \sim
: P\mapsto P $. It is mnemonic extension of the corresponding $q$
- definition by Kirschenhofer [24] and  Cigler [25]. This umbral
operator (without reference to to [24,25]) had been already used
in theoretical physics aiming at Quantum Mechanics on the lattice
[26]. The similar aim is represented by  [27,28] (see further
references there)  where incidence algebras are being prepared for
that purpose (Dirac notation included). As it is well known - the
classical umbral  [29,30]  and extended [10] finite operator
calculi  may be formulated in the reduced incidence algebra
language. Hence both applications of related tools to the same
goal are expected to meet at the arena of GHW algebra description
of both [14,13].

\end{document}